# Continuous $q$–Hermite polynomials:
# An elementary approach

*Johann Cigler*

Fakultät für Mathematik, Universität Wien
johann.cigler@univie.ac.at

**Abstract**

This overview article gives an elementary approach to continuous $q$–Hermite polynomials. We stress their relation to Fibonacci, Lucas and Chebyshev polynomials and to some $q$–analogues of these polynomials.

**0. Introduction**

Continuous $q$–Hermite polynomials are at the bottom of a large class of $q$–hypergeometric polynomials to which most of their properties can be generalized (cf. [2],[3],[11],[13]). This paper is the result of my efforts to understand the properties of these seminal polynomials and their relations to other polynomials. By "understand" I not only mean to verify complete proofs of all assertions but also to be aware of analogies with the classical case. Many mathematicians try to state and prove their theorems as general as possible. In doing so very often the roots of the theory or simpler and more illuminating proofs for special cases are disappearing from view. Therefore I want to give a direct approach to these polynomials without recourse to complicated $q$–hypergeometric identities or other more general theories. If there are other approaches in a similar manner please let me know.

It is interesting that there are two different kinds of Hermite polynomials in common use called respectively the physicists' and the probabilists' Hermite polynomials. This difference continues in the $q$–analogues. The analogues of the probabilistic Hermite polynomials are closely related to the Fibonacci and Lucas polynomials whereas the analogues of the physicists' polynomials can better be described using Chebyshev polynomials.

Some relations about these polynomials can easily be described using the notion of umbral inverse (cf.[16]). The sequence $(p_n(x))_{n \geq 0}$ of polynomials $p_n(x) = \sum_{k=0}^{n} a(n,k) x^k$ with $\deg p_n = n$ is umbrally inverse to the sequence $(r_n(x))_{n \geq 0}$ with $r_n(x) = \sum_{k=0}^{n} b(n,k) x^k$ if $\sum_{k=0}^{n} b(n,k) p_k(x) = x^n$. The umbral inverse of the Hermite polynomials $(H_n(x,s))$ is the sequence $(H_n(x,-s))$, but for the $q$–analogues the situation is much more complicated.

In order to make the paper self-contained I include a short introduction to $q$–identities.

The continuous $q$–Hermite polynomials $H_n(x | q)$ satisfy a $3$–term recurrence of Favard-type and are therefore orthogonal with respect to the linear functional $\Lambda_H$ on the polynomials defined by $\Lambda_H(H_n(x | q)) = [n = 0]$. (Here $[n = 0]$ denotes an Iverson bracket which is generally defined by $[P] = 1$ if property $P$ holds and $[P] = 0$ if $P$ does not hold). I use the method of Wm. Allaway [1] for determining the probability measure associated with the linear functional $\Lambda_H$.

Finally I introduce the Askey-Wilson $q$–differentiation operator and deduce the analogue of the Rodrigues formula.

Of course it is not possible to state all results from close to scratch. But I will at least sketch all relevant facts about the classical Hermite polynomials and about Fibonacci and Lucas polynomials or the equivalent Chebyshev polynomials. The reader who is not familiar with this background material can find proofs in the cited literature or other easily available sources.

## 1. Some background material about Hermite polynomials

### 1.1. Definitions

There are two closely related kinds of "Hermite polynomials" in common use.

The most popular ones are the polynomials $\mathbf{H}_n(x)$ which can be defined by the generating function
$$\sum_{n\geq 0}\frac{\mathbf{H}_n(x)}{n!}t^n = e^{2xt-t^2}. \tag{1.1}$$
They play an important role in applications and are often called the "physicists' Hermite polynomials".

The other kind are the "probabilists' Hermite polynomials" $H_n(x)$ with generating function
$$\sum_{n\geq 0}\frac{H_n(x)}{n!}t^n = e^{xt-\frac{t^2}{2}} \tag{1.2}$$
which are connected with the normal distribution in probability theory.

We consider more generally the bivariate Hermite polynomials $H_n(x,s)$ which can be defined by the generating function
$$\sum_{n\geq 0}\frac{H_n(x,s)}{n!}t^n = e^{tx-s\frac{t^2}{2}} \tag{1.3}$$

and contain both $H_n(x) := H_n(x,1)$ and $\mathbf{H}_n(x) = H_n(2x,2) = \sqrt{2^n}H_n\left(\sqrt{2}x,1\right)$ as special cases (cf. [2] or [16]).

Comparison of coefficients in
$$\sum_{n\geq 0}\frac{H_n(x,s)}{n!}t^n = e^{tx}e^{-s\frac{t^2}{2}} = \sum_j \frac{x^j t^j}{j!}\sum_k (-1)^k \frac{s^k t^{2k}}{2^k k!} = \sum_n \frac{t^n}{n!}\sum_{j+2k=n}\frac{n! x^j (-1)^k s^k}{j! 2^k k!}$$
gives

$$H_n(x,s) = \sum_{k=0}^{\left\lfloor\frac{n}{2}\right\rfloor}(-1)^k \binom{n}{2k}(2k-1)!! s^k x^{n-2k}. \tag{1.4}$$

The first polynomials are

$1,\ x,\ x^2 - s,\ x^3 - 3sx,\ x^4 - 6sx^2 + 3s^2,\ x^5 - 10sx^3 + 15s^2 x,\cdots.$

## 1.2. Some properties of the Hermite polynomials which will later be generalized

Since $\frac{d}{dt}e^{tx-s\frac{t^2}{2}} = (x-st)e^{tx-s\frac{t^2}{2}}$ we could alternatively define $H_n(x,s)$ as the coefficients of the uniquely determined formal power series solution $F(x,s,t) = \sum_{n\geq 0} \frac{H_n(x,s)}{n!}t^n$ of the equation

$$\frac{d}{dt}F(x,s,t) = (x-st)F(x,s,t) \tag{1.5}$$

with $F(x,s,0) = 1$.

This gives

$$H_{n+1}(x,s) = xH_n(x,s) - nsH_{n-1}(x,s) \tag{1.6}$$

with initial values $H_0(x,s) = 1$ and $H_1(x,s) = x$.

We also have

$$H_n(x+sD,s)1 = x^n. \tag{1.7}$$

The left-hand side should be interpreted in the following way: Replace in $H_n(x,s) = \sum_{k=0}^{n} h(n,k)x^k$ each power $x^k$ by the linear operator $(\mathbf{x}+sD)^k$ and apply the corresponding linear operator $\sum_{k=0}^{n} h(n,k)(\mathbf{x}+sD)^k$ to the constant polynomial 1. Here $D$ denotes the differentiation operator with respect to $x$ on the vector space of polynomials in $x$ and $\mathbf{x}$ the operator "multiplication by $x$". In order to simplify notation I shall denote the multiplication operator by $f$ on some vector space with the same symbol as the function $f$ itself.

Consider for example $H_3(x,s) = x^3 - 3sx$. We get

$$H_3(x+sD,s)1 = (x+sD)^3 1 - 3s(x+sD)1 = (x+sD)^2(x+sD)1 - 3s(x+sD)1$$
$$= (x+sD)(x+sD)x - 3sx = (x+sD)(x^2+s) - 3sx = x^3 + xs + 2sx - 3sx = x^3.$$

For the proof observe that (1.7) is obviously true for $n = 0$ and $n = 1$.

The general case follows by induction from

$$H_{n+1}(x+sD,s)1 = (x+sD)H_n(x+sD,s)1 - snH_{n-1}(x+sD,s)1 = (x+sD)x^n - nsx^{n-1} = x^{n+1}.$$

From (1.3) it is also clear that

$$\frac{d}{dx}H_n(x,s) = nH_{n-1}(x,s). \tag{1.8}$$

Thus (1.6) can also be written as

$$H_{n+1}(x,s) = (x - sD)H_n(x,s). \tag{1.9}$$

Therefore we have
$$H_n(x,s) = (x - sD)^n 1. \tag{1.10}$$

If we set
$$h_n(x,s) := (x + sD)^n 1 \tag{1.11}$$

then we get
$$h_n(x,s) = H_n(x,-s). \tag{1.12}$$

Equation (1.7) implies
$$\sum_{k=0}^{\lfloor \frac{n}{2} \rfloor} \binom{n}{2k}(2k-1)!! s^k H_{n-2k}(x,s) = x^n. \tag{1.13}$$

The operator $x - sD$ can also be written in the form
$$x - sD = e^{\frac{x^2}{2s}}(-sD)e^{-\frac{x^2}{2s}}. \tag{1.14}$$

For
$$e^{\frac{x^2}{2s}}(-sD)e^{-\frac{x^2}{2s}} f(x) = e^{\frac{x^2}{2s}}(-s)D\left(e^{-\frac{x^2}{2s}} f(x)\right) = e^{\frac{x^2}{2s}}(-s)\left(-\frac{x}{s}e^{-\frac{x^2}{2s}} + e^{-\frac{x^2}{2s}} f'(x)\right)$$
$$= x - sf'(x) = (x - sD)f(x).$$

This implies
$$(x - sD)^n = e^{\frac{x^2}{2s}}(-sD)^n e^{-\frac{x^2}{2s}} \tag{1.15}$$

and thus the so called **Rodrigues formula**
$$(-sD)^n e^{-\frac{x^2}{2s}} = H_n(x,s)e^{-\frac{x^2}{2s}}. \tag{1.16}$$

Let us also note some of the corresponding formulae for the polynomials $\mathbf{H}_n(x)$.

Here we have
$$D\mathbf{H}_n(x) = 2n\mathbf{H}_{n-1}(x), \tag{1.17}$$
$$\mathbf{H}_n(x) = (2x - D)^n 1, \tag{1.18}$$
$$\mathbf{H}_n(x) = (-1)^n e^{x^2} D^n e^{-x^2} \tag{1.19}$$

and
$$\mathbf{H}_n\left(x + \frac{1}{2}D\right) 1 = (2x)^n. \tag{1.20}$$

## 1.3. Inverse relations

Let $(A_n(x))_{n \geq 0}$ be a sequence of polynomials $A_n(x) = \sum_{k=0}^{n} u(n,k) x^k$ with $\deg A_n = n$ for all $n \in \mathbb{N}$.

Then there exist uniquely determined coefficients $v(n,k)$ such that $x^n = \sum_{k=0}^{n} v(n,k) A_k(x)$.

In this case we say that the sequence of polynomials $a_n(x) := \sum_{k=0}^{n} v(n,k) x^k$ is umbral inverse to $(A_n(x))_{n \geq 0}$. This is equivalent with $\left( v(i,j) \right)_{i,j=0}^{n} = \left( \left( u(i,j) \right)_{i,j=0}^{n} \right)^{-1}$.

Umbral inverse polynomials give rise to inverse relations: Let $(x(n))_{n \geq 0}$ be any sequence of numbers such that $x(0) = 1$ and let $y(n) = \sum_{k=0}^{n} u(n,k) x(k)$. Then $x(n) = \sum_{k=0}^{n} v(n,k) y(k)$ and vice versa. To see this we have only to apply the linear functional $\Lambda$ defined by $\Lambda(x^n) = x(n)$.

Equation (1.13) says that the sequence $h_n(x,s) = H_n(x,-s)$ is umbral inverse to $H_n(x,s)$.

## 1.4. Moments and orthogonality

If we define a linear functional $\Lambda_H$ on the vector space of polynomials in $x$ by

$$\Lambda_H(H_n(x,s)) = [n = 0] \qquad (1.21)$$

then (1.13) implies that $\Lambda_H(x^{2n+1}) = 0$ and

$$\Lambda_H(x^{2n}) = (2n-1)!! s^n. \qquad (1.22)$$

The numbers $\Lambda_H(x^n)$ are the so called **moments** of the Hermite polynomials.

The Hermite polynomials are **orthogonal polynomials** on the real line.

A sequence of monic polynomials $p_n(x)$ with degree $\deg p_n = n$ is called orthogonal with respect to a linear functional $\Lambda$ if $\Lambda\left( p_m(x) \overline{p_n(x)} \right) = C_n [m = n]$ with $C_n \neq 0$.

Since the sequence $(p_n(x))_{n \geq 0}$ is a basis for the vector space of polynomials the linear functional $\Lambda$ is uniquely determined by $\Lambda(p_n(x)) = [n = 0]$.

By Favard's theorem (cf. e.g. [2]) orthogonality on the real line is equivalent with a $3-$term recurrence

$$p_{n+1}(x) = (x + a_n) p_n(x) + b_n p_{n-1}(x) \qquad (1.23)$$

with $p_{-1}(x) = 0$ and $b_n \neq 0$.

Let us note that for $H_n(x, s)$ with $s > 0$ the linear functional $\Lambda_H$ is given by a probability measure. More precisely we have

$$\Lambda_H(f) = \frac{1}{\sqrt{2\pi s}} \int_{-\infty}^{\infty} f(x) e^{-\frac{x^2}{2s}} dx. \tag{1.24}$$

To see this let us start from the well-known Fourier transform of the Gaussian measure

$$\frac{1}{\sqrt{2\pi s}} \int_{-\infty}^{\infty} e^{iyx} e^{-\frac{x^2}{2s}} dx = e^{-s\frac{y^2}{2}}. \tag{1.25}$$

If we differentiate $2n$ times with respect to $y$ and let afterwards $y \to 0$ then we get by (1.16) with $\frac{1}{s}$ in place of $s$

$$\frac{1}{\sqrt{2\pi s}} \int_{-\infty}^{\infty} (ix)^{2n} e^{-\frac{x^2}{2s}} dx = \left(\frac{d}{dy}\right)^{2n} e^{-s\frac{y^2}{2}} \bigg|_{y=0} = (-s)^{2n} H_{2n}(0, \frac{1}{s}) = (-1)^n s^n (2n-1)!!$$

or

$$\frac{1}{\sqrt{2\pi s}} \int_{-\infty}^{\infty} x^{2n} e^{-\frac{x^2}{2s}} dx = s^n (2n-1)!!.$$

It is clear that $\frac{1}{\sqrt{2\pi s}} \int_{-\infty}^{\infty} x^{2n+1} e^{-\frac{x^2}{2s}} dx = 0.$

Comparing with (1.22) we get

$$\Lambda_H(x^n) = \frac{1}{\sqrt{2\pi s}} \int_{-\infty}^{\infty} x^n e^{-\frac{x^2}{2s}} dx, \tag{1.26}$$

which implies (1.24) for polynomials $f(x)$. The orthogonality relation is

$$\Lambda_H(H_m(x,s) H_n(x,s)) = \frac{1}{\sqrt{2\pi s}} \int_{-\infty}^{\infty} H_m(x,s) H_n(x,s) e^{-\frac{x^2}{2s}} dx = s^n n! [m=n]. \tag{1.27}$$

Since $\Lambda_H(H_m(x,s) H_n(x,s)) = \Lambda_H(x^m H_n(x,s))$ for $m \leq n$ we see from the recurrence relation (1.6) that

$$\Lambda_H(x^n H_n(x)) = \Lambda_H(x^{n-1} H_{n+1}(x)) + ns \Lambda_H(x^{n-1} H_n(x)) = ns \Lambda_H(x^{n-1} H_n(x))$$

and by iteration $\Lambda_H(x^n H_n(x)) = n! s^n$.

For the physicists' Hermite polynomials we get

$$\frac{1}{\sqrt{\pi}} \int_{-\infty}^{\infty} \mathbf{H}_m(x) \mathbf{H}_n(x) e^{-x^2} dx = 2^n n! [m=n]. \tag{1.28}$$

## 2. Some definitions and results from q-analysis

### 2.1. Preliminaries

Let $q \neq 1$ be a real number and let $[n]_q = \dfrac{1-q^n}{1-q}$, $[n]_q! = [1]_q[2]_q \cdots [n]_q$ with $[0]_q! = 1$

and $\begin{bmatrix} n \\ k \end{bmatrix}_q = \dfrac{[n]_q!}{[k]_q![n-k]_q!}$ for $0 \leq k \leq n$ and $\begin{bmatrix} n \\ k \end{bmatrix}_q = 0$ else.

Mostly we will drop the index $q$ and simply write $[n]$ instead of $[n]_q$.

For $q \to 1$ we have $[n]_q \to n$ and $\begin{bmatrix} n \\ k \end{bmatrix}_q \to \binom{n}{k}$. Therefore $[n]_q$ is called a $q$-analogue of $n$ and $\begin{bmatrix} n \\ k \end{bmatrix}_q$ a $q$-analogue of $\binom{n}{k}$.

The $q$-binomial coefficients satisfy

$$\begin{bmatrix} n+1 \\ k \end{bmatrix} = q^k \begin{bmatrix} n \\ k \end{bmatrix} + \begin{bmatrix} n \\ k-1 \end{bmatrix} \tag{2.1}$$

and

$$\begin{bmatrix} n+1 \\ k \end{bmatrix} = \begin{bmatrix} n \\ k \end{bmatrix} + q^{n+1-k} \begin{bmatrix} n \\ k-1 \end{bmatrix}. \tag{2.2}$$

We will also use the notation $(x;q)_n = (1-x)(1-qx)\cdots(1-q^{n-1}x)$ and
$(x;q)_\infty = (1-x)(1-qx)(1-q^2x)\cdots$.

Let $D_q$ be the $q$-differentiation operator on the vector space of polynomials in $x$ defined by
$D_q p(x) = \dfrac{p(x) - p(qx)}{(1-q)x}$. It is uniquely determined by $D_q x^n = [n]_q x^{n-1}$ for all $n \in \mathbb{N}$. For $q \to 1$ it converges to ordinary differentiation.

Note that
$$D_q(f(x)g(x)) = g(x)D_q(f(x)) + f(qx)D_q(g(x)). \tag{2.3}$$

For this and the following sections the interested reader could also consult [2], [5], [6] or [11].

### 2.2. Some simple q-analogues of the binomial theorem

Let $\varepsilon = \varepsilon_q$ be the linear operator on the polynomials defined by $\varepsilon_q p(x) = p(qx)$.
Then the following $q$-analogue of the binomial theorem holds:

$$\left(x + \varepsilon_q\right)^n = \sum_{k=0}^{n} \begin{bmatrix} n \\ k \end{bmatrix}_q x^k \varepsilon_q^{n-k}. \tag{2.4}$$

To prove this observe that $\varepsilon_q x = qx\varepsilon_q$ because $\varepsilon_q xp(x) = qxp(qx) = qx\varepsilon_q p(x)$ for all polynomials $p(x)$.

Let more generally $A$ and $B$ be linear operators on the polynomials satisfying

$$BA = qAB \tag{2.5}$$

then

$$(A+B)^n = \sum_{k=0}^{n} \begin{bmatrix} n \\ k \end{bmatrix} A^k B^{n-k}. \tag{2.6}$$

This can easily be proved by induction:
$$(A+B)\sum_{k=0}^{n} \begin{bmatrix} n \\ k \end{bmatrix} A^k B^{n-k} = \sum_{k=0}^{n} \begin{bmatrix} n \\ k \end{bmatrix} A^{k+1} B^{n-k} + \sum_{k=0}^{n} \begin{bmatrix} n \\ k \end{bmatrix} q^k A^k B^{n-k+1}$$
$$= \sum_{k=0}^{n+1} \begin{bmatrix} n \\ k-1 \end{bmatrix} A^k B^{n+1-k} + \sum_{k=0}^{n} \begin{bmatrix} n \\ k \end{bmatrix} q^k A^k B^{n-k+1} = \sum_{k=0}^{n+1} \left(\begin{bmatrix} n \\ k-1 \end{bmatrix} + q^k \begin{bmatrix} n \\ k \end{bmatrix}\right) A^k B^{n-k+1} = \sum_{k=0}^{n+1} \begin{bmatrix} n+1 \\ k \end{bmatrix} A^k B^{n+1-k}.$$

If we apply the operator $(x+y\varepsilon_q)^n$ to the constant polynomial $p(x)=1$ we get the **Rogers-Szegö polynomials**

$$R_n(x,y,q) = (x+y\varepsilon_q)^n 1 = \sum_{k=0}^{n} \begin{bmatrix} n \\ k \end{bmatrix}_q x^k y^{n-k}. \tag{2.7}$$

They satisfy

$$D_q R_n(x,y,q) = [n]_q R_{n-1}(x,y,q) \tag{2.8}$$

and

$$R_n(x,y,q) = (x+y)R_{n-1}(x,y,q) + (q^{n-1}-1)xyR_{n-2}(x,y,q). \tag{2.9}$$

For

$$D_q R_n(x,y,q) = \sum_{k=0}^{n} \begin{bmatrix} n \\ k \end{bmatrix} [k] x^{k-1} y^{n-k} = [n]\sum_{k=0}^{n} \begin{bmatrix} n-1 \\ k-1 \end{bmatrix} x^{k-1} y^{n-k} = [n]\sum_{k=0}^{n-1} \begin{bmatrix} n-1 \\ k \end{bmatrix} x^k y^{n-1-k} = [n]R_{n-1}(x,y,q)$$

and

$$\varepsilon_q = 1 + (q-1)xD_q. \tag{2.10}$$

The last identity holds because
$(1+(q-1)xD_q)x^n = x^n + (q-1)x[n]x^{n-1} = x^n(1+q^n-1) = (qx)^n = \varepsilon_q x^n.$
This finally gives

$$R_n(x,y,q) = (x+y(1+(q-1)xD_q))R_{n-1}(x,y,q) = (x+y)R_{n-1}(x,y,q) + (q^{n-1}-1)xyR_{n-2}(x,y,q).$$

It should be noted that recurrences such as relation (2.9) can also be obtained from the $q$-Zeilberger algorithm. The implementation by Peter Paule and Axel Riese [14] gives rise to a computer proof:

```
qZeil[qBinomial[n, k, q] x^k y^(n - k), {k, 0, n}, n, 2]
SUM[n] == - (1 - q^(-1+n)) x y SUM[-2 + n] + (x + y) SUM[-1 + n]
```

There are some special cases known where $R_n(x, y, q)$ has simple values.

The first one is the Gauss formula $R_{2n}(1,-1,q) = (1-q)(1-q^3)\cdots(1-q^{2n-1})$ and $R_{2n+1}(1,-1,q) = 0$.

The second one is

$$R_n(q,1,q^2) = \sum_{k=0}^{n} \begin{bmatrix} n \\ k \end{bmatrix}_{q^2} q^k = (1+q)(1+q^2)\cdots(1+q^n).  \tag{2.11}$$

A third one is

$$R_n(-q,1,q) = \sum_{k=0}^{n} \begin{bmatrix} n \\ k \end{bmatrix}(-q)^k = (1-q)\left(1-q^3\right)\cdots\left(1-q^{2\left\lfloor\frac{n+1}{2}\right\rfloor-1}\right).  \tag{2.12}$$

These results can easily be verified by using the recurrence (2.9).

The simplest $q$-analogue of the binomial theorem can be obtained by setting $A = x\varepsilon_q$ and $B = y\varepsilon_q$.

This gives

$$(y+x)(y+qx)\cdots(y+q^{n-1}x) = \sum_{k=0}^{n} q^{\binom{k}{2}} \begin{bmatrix} n \\ k \end{bmatrix} x^k y^{n-k}.  \tag{2.13}$$

For

$$(x\varepsilon + y\varepsilon)^n = (x+y)\varepsilon(x+y)\varepsilon\cdots(x+y)\varepsilon = (x+y)(qx+y)(q^2x+y)\cdots(q^{n-1}x+y)\varepsilon^n.$$

The computer proof gives

```
Simplify[qZeil[q^Binomial[k, 2] qBinomial[n, k, q] x^k y^(n - k), {k, 0, n}, n, 1]]
(q^(-1+n) x + y) SUM[-1 + n] == SUM[n]
```

## 2.3. The q-exponential series

Let

$$e(t) = e_q(t) = \sum_{n} \frac{t^n}{[n]_q!}  \tag{2.14}$$

denote the $q$-exponential series.

For $q = 0$ we get $e_0(t) = \sum_{n \geq 0} t^n = \frac{1}{1-t}$. Thus $e_q(t)$ for $0 \leq t \leq 1$ is some sort of interpolation between the geometric series and the usual exponential series.

We can consider (2.14) as a formal power series. It is clear that $D_q e_q(ax) = a e_q(ax)$ if as usual for formal power series we apply $D_q$ to each term.

An important property is
$$e_q(At)e_q(Bt) = e_q((A+B)t) \tag{2.15}$$
if $BA = qAB$.

This follows from
$$e_q(At)e_q(Bt) = \sum_j \frac{A^j t^j}{[j]!} \sum_k \frac{B^k t^k}{[k]!} = \sum_n \frac{t^n}{[n]!} \sum_{j=0}^n \begin{bmatrix} n \\ j \end{bmatrix} A^j B^{n-j} = \sum_n \frac{t^n}{[n]!}(A+B)^n = e_q((A+B)t).$$

For $A = x$ and $B = y\varepsilon_q$ (2.15) applied to the constant polynomial 1 gives
$$\sum_n \frac{R_n(x,y,q)}{[n]!} t^n = e(xt)e(yt). \tag{2.16}$$

If $A = x$ and $B = -x\varepsilon_q$ then we have also $BA = qAB$.
Therefore
$$e_q(xt)e_q(-x\varepsilon_q t) = e_q\left(x(1-\varepsilon_q)t\right).$$

If we apply this to the constant polynomial 1 we get $e(xt)\sum_n (-1)^n q^{\binom{n}{2}} \frac{x^n t^n}{[n]!} = 1$ or equivalently
$$\frac{1}{e(t)} = \sum_n (-1)^n q^{\binom{n}{2}} \frac{t^n}{[n]!}. \tag{2.17}$$

Since $e(xt) = \frac{1}{e(yt)} \sum_n \frac{R_n(x,y,q)}{[n]!} t^n$ we get by comparing coefficients
$$\sum_{k=0}^n (-y)^{n-k} q^{\binom{n-k}{2}} \begin{bmatrix} n \\ k \end{bmatrix} R_k(x,y,q) = x^n. \tag{2.18}$$

Therefore the umbral inverse sequence to the Rogers-Szegö polynomials is the sequence
$$r_n(x,y,q) = \sum_{k=0}^n (-y)^{n-k} q^{\binom{n-k}{2}} \begin{bmatrix} n \\ k \end{bmatrix} x^k = (x-y)(x-qy)\cdots(x-q^{n-1}y). \tag{2.19}$$

As a consequence the moments of the Rogers-Szegö polynomials are
$$\Lambda_{R(x,y,q)}(x^n) = (-y)^n q^{\binom{n}{2}}. \tag{2.20}$$

## 2.4. Rogers-Szegö polynomials on the unit circle of the complex plane

Since the recurrence (2.9) of the Rogers-Szegö polynomials is not of the form (1.23) they are not orthogonal on the real line with respect to $\Lambda_{R(x,y,q)}$.

But consider the functions $R_n(z, y, q) = R_n(e^{i\vartheta}, y, q)$ on the unit circle $z = e^{i\vartheta}$ with $0 \leq \vartheta < 2\pi$. These are trigonometric polynomials, i.e. finite linear combinations of $e^{ik\vartheta}$ for $k \in \mathbb{Z}$.

We want to find some $y$ such that the polynomials $R_n(z, y, q)$ are orthogonal with respect to a probability measure on the unit circle. The simplest example is (cf. [17])

$$R_n(z) = R_n(z, -\sqrt{q}, q) = \sum_{k=0}^{n} (-1)^{n-k} \begin{bmatrix} n \\ k \end{bmatrix} q^{\frac{n-k}{2}} z^k \tag{2.21}$$

for $0 < q < 1$.

Consider the function $v(\vartheta, q) = \frac{1}{\sqrt{2\pi s}} \sum_{j=-\infty}^{\infty} e^{-\frac{(\vartheta - 2\pi j)^2}{2s}}$. This is the density of the so called wrapped Gaussian measure on the unit circle. For $\varphi(z) = \varphi(e^{i\vartheta})$ we let $\varphi(e^{ix})$ be periodic with period $2\pi$ on $\mathbb{R}$. Then we get

$$\frac{1}{\sqrt{2\pi s}} \int_{-\infty}^{\infty} \varphi(e^{ix}) e^{-\frac{x^2}{2s}} dx = \frac{1}{\sqrt{2\pi s}} \sum_{j=-\infty}^{\infty} \int_{0}^{2\pi} \varphi(e^{i\vartheta}) e^{-\frac{(\vartheta - 2\pi j)^2}{2s}} d\vartheta = \int_{0}^{2\pi} \varphi(e^{i\vartheta}) v(\vartheta, q) d\vartheta.$$

For $\varphi(z) = z^n$, $n \in \mathbb{Z}$, we get by (1.25) $\frac{1}{\sqrt{2\pi s}} \int_{-\infty}^{\infty} e^{inx} e^{-\frac{x^2}{2s}} dx = e^{-s\frac{n^2}{2}}$.

Choose now $s$ such that $q = e^{-s}$.
Then

$$\Lambda_R(z^n) = \int_{0}^{2\pi} e^{in\vartheta} v(\vartheta, q) d\vartheta = q^{\frac{n^2}{2}}. \tag{2.22}$$

The polynomials $R_n(z, -\sqrt{q}, q)$ on the unit circle are orthogonal with respect to $\Lambda_R$.
More precisely

$$\Lambda_R\left(R_m(z) \overline{R_n(z)}\right) = (q; q)_n [n = m]. \tag{2.23}$$

Let $I_{m,n} = \Lambda_R\left(R_m(z) \overline{R_n(z)}\right)$.

Then

$$I_{m,n} = (-q)^{\frac{m+n}{2}} \sum_{j=0}^{m} (-1)^j q^{\binom{j}{2}} \begin{bmatrix} m \\ j \end{bmatrix} \prod_{\ell=0}^{n-1} \left(1 - q^{\ell - j}\right). \tag{2.24}$$

For

$$I_{m,n} = \Lambda_R\left(R_m(z)\overline{R_n(z)}\right) = \Lambda_R\left(R_m(z)R_n(\overline{z})\right) = \Lambda_R\left(\sum_{j=0}^{m}(-1)^{m-j}\begin{bmatrix}m\\j\end{bmatrix}q^{\frac{m-j}{2}}z^j \sum_{k=0}^{n}(-1)^{n-k}\begin{bmatrix}n\\k\end{bmatrix}q^{\frac{n-k}{2}}z^{-k}\right)$$

$$= \sum_{j=0}^{m}\sum_{k=0}^{n}(-1)^{m+n-j-k}\begin{bmatrix}m\\j\end{bmatrix}\begin{bmatrix}n\\k\end{bmatrix}q^{\frac{m-j}{2}}q^{\frac{n-k}{2}}\Lambda_R(z^{j-k}) = \sum_{j=0}^{m}\sum_{k=0}^{n}(-1)^{m+n-j-k}\begin{bmatrix}m\\j\end{bmatrix}\begin{bmatrix}n\\k\end{bmatrix}q^{\frac{m-j}{2}}q^{\frac{n-k}{2}}q^{\frac{(j-k)^2}{2}}$$

$$= (-1)^{m+n}q^{\frac{m+n}{2}}\sum_{j=0}^{m}(-1)^j q^{\binom{j}{2}}\begin{bmatrix}m\\j\end{bmatrix}\sum_{k=0}^{n}(-1)^k q^{\binom{k}{2}}\begin{bmatrix}n\\k\end{bmatrix}(q^{-j})^k = (-1)^{m+n}q^{\frac{m+n}{2}}\sum_{j=0}^{m}(-1)^j q^{\binom{j}{2}}\begin{bmatrix}m\\j\end{bmatrix}\prod_{\ell=0}^{n-1}(1-q^{\ell-j}).$$

For $m < n$ we have $\prod_{\ell=0}^{n-1}(1-q^{\ell-j}) = 0$ for each $j$ with $0 \le j \le m$.

For $n = m$ we get

$$I_{n,n} = q^n \sum_{j=0}^{n}(-1)^j q^{\binom{j}{2}}\begin{bmatrix}n\\j\end{bmatrix}\prod_{\ell=0}^{n-1}(1-q^{\ell-j}) = q^n(-1)^n q^{\binom{n}{2}}\prod_{\ell=0}^{n-1}(1-q^{\ell-n}) = (q;q)_n.$$

Note that this result has no counterpart for $q = 1$.

### 2.5. Jacobi's triple product identity

Instead of the formal power series $e_q(t)$ we can for $|q| < 1$ also consider the convergent series $e_q\left(\frac{t}{1-q}\right) = \sum_n \frac{t^n}{(q;q)_n}$ which converges for $|t| < 1$. This can also be considered as a sort of $q$-analogue of the exponential function although there is no limit for $q \to 1$.

It satisfies

$$\sum_n \frac{t^n}{(q;q)_n} = \frac{1}{(t;q)_\infty}. \tag{2.25}$$

For $(1-t)\sum_n \frac{t^n}{(q;q)_n} = \sum_n \frac{t^n}{(q;q)_n} - \sum_n \frac{t^{n+1}}{(q;q)_n} = \sum_n \frac{t^n}{(q;q)_n} - \sum_n \frac{t^n(1-q^n)}{(q;q)_n} = \sum_n \frac{(qt)^n}{(q;q)_n}$

and by iterating

$(t;q)_k \sum_n \frac{t^n}{(q;q)_n} = \sum_n \frac{(q^k t)^n}{(q;q)_n}.$

For $k \to \infty$ the right-hand side converges to 1.

In the same way we get

$$\sum_n \frac{(-1)^n q^{\binom{n}{2}} t^n}{(q;q)_n} = (t;q)_\infty. \tag{2.26}$$

This could also be obtained from (2.13):

If we let $n \to \infty$ in

$$(t;q)_n = (1-t)(1-qt)\cdots(1-q^{n-1}t) = \sum_{k=0}^{n}(-1)^k q^{\binom{k}{2}} \begin{bmatrix} n \\ k \end{bmatrix} t^k$$

the left-hand side converges to $(t;q)_\infty$ and the right-hand side to $\sum_{k}\frac{(-1)^k q^{\binom{k}{2}} t^k}{(q;q)_k}$ because

$$\lim_{k \to \infty}\begin{bmatrix} n \\ k \end{bmatrix}_q = \lim_{k \to \infty}\frac{(q;q)_n}{(q;q)_k (q;q)_{n-k}} = \frac{(q;q)_\infty}{(q;q)_k (q;q)_\infty} = \frac{1}{(q;q)_k}.$$

A similar argument leads to the **triple product identity of Jacobi**

$$\sum_{k \in \mathbb{Z}}(-1)^k q^{\binom{k}{2}} x^k = (x;q)_\infty \left(\frac{q}{x};q\right)_\infty (q;q)_\infty. \tag{2.27}$$

Consider

$$(1-q^{-n}x)(1-q^{-n+1}x)\cdots(1-q^{-1}x)(1-x)(1-qx)\cdots(1-q^{n-1}x) = (q^{-n}x;q)_{2n}$$

$$= \sum_{k=0}^{2n}(-1)^k q^{\binom{k}{2}}\begin{bmatrix} 2n \\ k \end{bmatrix}(q^{-n}x)^k = \sum_{j=-n}^{n}(-1)^{n+j} q^{\binom{n+j}{2}-n(j+n)}\begin{bmatrix} 2n \\ n+j \end{bmatrix} x^{n+j}.$$

Multiplying both sides with $q^{\binom{n+1}{2}}x^{-n} = \frac{q^n}{x}\frac{q^{n-1}}{x}\cdots\frac{q}{x}$ and observing that

$\binom{n+j}{2} - n(n+j) + \binom{n+1}{2} = \binom{j}{2}$ we get

$$\left(1-\frac{q^n}{x}\right)\left(1-\frac{q^{n-1}}{x}\right)\cdots\left(1-\frac{q}{x}\right)(1-x)(1-qx)\cdots(1-q^{n-1}x) = \sum_{j=-n}^{n}(-1)^j q^{\binom{j}{2}}\begin{bmatrix} 2n \\ n+j \end{bmatrix} x^j.$$

For $n \to \infty$ we obtain (2.27).

More details about these themes can be found in [2], [5], [6] and [11].

## 3. Fibonacci and Lucas polynomials

### 3.1. Fibonacci polynomials

Let us state some well-known facts about these polynomials. Details may be found in [8].
The **Fibonacci polynomials**

$$F_{n+1}(x,s) = \sum_{k=0}^{\lfloor \frac{n}{2} \rfloor} \binom{n-k}{k} s^k x^{n-2k} \tag{3.1}$$

satisfy

$F_{n+1}(x,s) = xF_n(x,s) + sF_{n-2}(x,s)$ with initial values $F_0(x,s) = 0$ and $F_1(x,s) = 1$ and Binet's formula

$$F_n(x,s) = \frac{u^n - v^n}{u - v} \text{ with } u = \frac{x + \sqrt{x^2 + 4s}}{2} \text{ and } v = \frac{x - \sqrt{x^2 + 4s}}{2}.$$

For $(x,s) = (1,1)$ we get the Fibonacci numbers. Their most important property is $\gcd(F_m, F_n) = F_{\gcd(m,n)}$ if $\gcd(a,b)$ denote the greatest common divisor of the numbers $a, b$.

In our context some people would prefer to use instead the polynomials $f_n(x,s) = F_{n+1}(x,s)$ with $\deg f_n = n$. But this would clash with the above result for Fibonacci numbers.

### 3.2. Lucas polynomials

The **Lucas polynomials**

$$L_n(x,s) = \sum_{k=0}^{\lfloor \frac{n}{2} \rfloor} \frac{n}{n-k} \binom{n-k}{k} s^k x^{n-2k} \tag{3.2}$$

satisfy
$L_{n+1}(x,s) = xL_n(x,s) + sL_{n-2}(x,s)$ with initial values $L_0(x,s) = 2$ and $L_1(x,s) = x$ and Binet's formula

$$L_n(x,s) = u^n + v^n \text{ with } u = \frac{x + \sqrt{x^2 + 4s}}{2} \text{ and } v = \frac{x - \sqrt{x^2 + 4s}}{2}.$$

These polynomials are related by
$$L_n(x,s) = F_{n+1}(x,s) + sF_{n-1}(x,s) \tag{3.3}$$
for $n > 0$.

In order to simply formulate (3.5) we denote by $L_n^*(x,s)$ the polynomials which coincide with $L_n(x,s)$ for $n > 0$ but satisfy $L_0^*(x,s) = 1$.

Since the polynomials $F_{n+1}(x,s)$ and $L_n^*(x,s)$ have degree $n$ it is clear that they form a basis for the vector space of polynomials. Therefore $x^n$ can be expressed as a linear combination of these polynomials. More precisely we have

$$\sum_{k=0}^{\lfloor \frac{n}{2} \rfloor} \left( \binom{n}{k} - \binom{n}{k-1} \right) s^k F_{n+1-2k}(x,-s) = x^n \tag{3.4}$$

and

$$\sum_{k=0}^{\lfloor \frac{n}{2} \rfloor} \binom{n}{k} s^k L_{n-2k}^*(x,-s) = x^n. \tag{3.5}$$

Let us first prove (3.5). Let $\alpha = \dfrac{x + \sqrt{x^2 - 4s}}{2}$ and $\beta = \dfrac{x - \sqrt{x^2 - 4s}}{2}$.

For odd $n$ we have

$$x^n = (\alpha + \beta)^n = \sum_{k=0}^{n} \binom{n}{k} \alpha^k \beta^{n-k} = \sum_{k=0}^{\lfloor n/2 \rfloor} \binom{n}{k} \left( \alpha^k \beta^{n-k} + \alpha^{n-k} \beta^k \right) = \sum_{k=0}^{\lfloor n/2 \rfloor} \binom{n}{k} \alpha^k \beta^k \left( \beta^{n-2k} + \alpha^{n-2k} \right)$$

$$= \sum_{k=0}^{\lfloor n/2 \rfloor} \binom{n}{k} s^k \left( \beta^{n-2k} + \alpha^{n-2k} \right) = \sum_{k=0}^{\lfloor n/2 \rfloor} \binom{n}{k} s^k L_{n-2k}(x, s).$$

For $n = 2m$ we get

$$x^{2m} = \sum_{k=0}^{2m} \binom{2m}{k} \alpha^k \beta^{2m-k} = \binom{2m}{m} \alpha^m \beta^m + \sum_{k=0}^{m-1} \binom{2m}{k} \left( \alpha^k \beta^{2m-k} + \alpha^{2m-k} \beta^k \right)$$

$$= \sum_{k=0}^{m-1} \binom{2m}{k} s^k \left( \beta^{2m-2k} + \alpha^{2m-2k} \right) + \binom{2m}{m} s^m = \sum_{k=0}^{m} \binom{2m}{k} s^k L^*_{2m-2k}(x, s).$$

Using (3.3) we see that (3.5) implies (3.4).

These identities immediately imply some well-known inverse relations:

Let $(a(n))_{n \geq 0}$ and $(b(n))_{n \geq 0}$ be two sequences with $a(0) = b(0) = 1$ which satisfy

$$b(n) = \sum_{k=0}^{\lfloor n/2 \rfloor} \binom{n-k}{k} s^k a(n-2k). \tag{3.6}$$

Then

$$a(n) = \sum_{k=0}^{\lfloor n/2 \rfloor} \left( \binom{n}{k} - \binom{n}{k-1} \right) (-s)^k b(n-2k). \tag{3.7}$$

If

$$b(n) = \sum_{k=0}^{\lfloor n/2 \rfloor} \frac{n}{n-k} \binom{n-k}{k} s^k a(n-2k). \tag{3.8}$$

Then

$$a(n) = \sum_{k=0}^{\lfloor n/2 \rfloor} \binom{n}{k} (-s)^k b(n-2k). \tag{3.9}$$

These identities allow us also to compute the moments corresponding to the Fibonacci and Lucas polynomials which turn out to be Catalan numbers and central binomial coefficients.

Let $\Lambda_F$ be defined by $\Lambda_F(F_{n+1}(x,-1)) = [n = 0]$. Then (3.4) implies

$$\Lambda_F(x^{2n}) = \frac{1}{n+1}\binom{2n}{n} = C_n. \tag{3.10}$$

From this it is easy to verify the formula

$$\Lambda_F(x^n) = \frac{1}{2\pi}\int_{-2}^{2} x^n \sqrt{4-x^2}\,dx. \tag{3.11}$$

If we define $\Lambda_L$ by $\Lambda_L(L_n(x,-1)) = [n = 0]$ then (3.5) implies

$$\Lambda_L(x^{2n}) = \binom{2n}{n}. \tag{3.12}$$

In this case

$$\Lambda_L(x^n) = \frac{1}{\pi}\int_{-2}^{2} \frac{x^n}{\sqrt{4-x^2}}\,dx. \tag{3.13}$$

## 4. Chebyshev polynomials

Chebyshev polynomials are Fibonacci and Lucas polynomials in disguised form. They are often preferred to them because of their close relations to the trigonometric functions. More details may be found in [9].

The **Chebyshev polynomials of the first kind** $T_n(x)$ are a variant of the Lucas polynomials

$$T_n(x) = \frac{1}{2}L_n(2x,-1) = 2^{n-1}L_n\left(x,-\frac{1}{4}\right) = \frac{1}{2}\left(\left(x+\sqrt{x^2-1}\right)^n + \left(x-\sqrt{x^2-1}\right)^n\right). \tag{4.1}$$

They are characterized by

$$T_n(\cos\vartheta) = \cos n\vartheta. \tag{4.2}$$

The **Chebyshev polynomials of the second kind** $U_n(x)$ are a variant of the Fibonacci polynomials

$$U_n(x) = F_{n+1}(2x,-1) = 2^n F_{n+1}\left(x,-\frac{1}{4}\right). \tag{4.3}$$

They are characterized by

$$U_n(\cos\vartheta) = \frac{\sin(n+1)\vartheta}{\sin\vartheta}. \tag{4.4}$$

Both types are connected by the identity

$$\left(x+\sqrt{x^2-1}\right)^n = T_n(x) + U_{n-1}(x)\sqrt{x^2-1} \tag{4.5}$$

if we set $U_{-1}(x) = 0$.

This identity is nothing else than $e^{in\vartheta} = (\cos\vartheta + i\sin\vartheta)^n$ in disguised form.

The linear functional $\Lambda_T$ corresponding to the polynomials of the first kind satisfies

$$\Lambda_T(f(x)) = \frac{1}{\pi}\int_{-1}^{1}\frac{f(x)}{\sqrt{1-x^2}}dx. \tag{4.6}$$

The moments are $\Lambda_T(x^{2n+1}) = 0$ and

$$\Lambda_T(x^{2n}) = \frac{1}{2^{2n}}\binom{2n}{n}. \tag{4.7}$$

The linear functional $\Lambda_U$ corresponding to the polynomials of the second kind satisfies

$$\Lambda_U(f(x)) = \frac{2}{\pi}\int_{-1}^{1} f(x)\sqrt{1-x^2}\,dx. \tag{4.8}$$

The moments are $\Lambda_U(x^{2n+1}) = 0$ and

$$\Lambda_U(x^{2n}) = \frac{1}{2^{2n}}\frac{1}{(n+1)}\binom{2n}{n}. \tag{4.9}$$

The inverse relations are

$$\frac{1}{2^n}\sum_{k=0}^{\lfloor\frac{n}{2}\rfloor}\left(\binom{n}{k}-\binom{n}{k-1}\right)U_{n-2k}(x) = x^n \tag{4.10}$$

and

$$\frac{1}{2^{n-1}}\sum_{k=0}^{\lfloor\frac{n}{2}\rfloor}\binom{n}{k}T^*_{n-2k}(x) = x^n \tag{4.11}$$

if $T_0^*(x) = \frac{1}{2}$ and $T_n^*(x) = T_n(x)$ for $n > 0$.

## 5. Generalities about q-Hermite polynomials

### 5.1. Definition

There are many $q$-analogues of the Hermite polynomials with interesting properties. We shall consider some variants of the so called **continuous $q$-Hermite polynomials** $H_n(x\,|\,q)$. Their main properties can be found in [1], [2], [3], [10],[11],[12],[13] and [20].

Let $\dfrac{\partial}{\partial t} F(t) = \dfrac{F(qt)-F(t)}{(q-1)t}$ denote $q$-differentiation with respect to the variable $t$ for formal power series.

In view of (1.5) the most natural $q$-analogue of $H_n(x,s)$ is $\tilde{H}_n(x,s,q)$ which satisfies

$$\frac{\partial}{\partial t} \sum_{n \geq 0} \frac{\tilde{H}_n(x,s,q)}{[n]!} t^n = (x - st) \sum_{n \geq 0} \frac{\tilde{H}_n(x,s,q)}{[n]!} t^n \tag{5.1}$$

or equivalently

$$\tilde{H}_{n+1}(x,s,q) = x\tilde{H}_n(x,s,q) - [n]s\tilde{H}_{n-1}(x,s,q) \tag{5.2}$$

which for $q \to 1$ converges to $H_n(x,s)$.

This $q$-analogue is important in combinatorics (cf. [12]) but for our purposes the closely related polynomials $H_n(x,s,q) = \tilde{H}_n(x,(1-q)s,q)$ are better suited.

Therefore we define the **bivariate $q$-Hermite polynomials** $H_n(x,s,q)$ by

$$H_{n+1}(x,s,q) = xH_n(x,s,q) - \left(1-q^n\right)sH_{n-1}(x,s,q) \tag{5.3}$$

with initial values $H_0(x,s,q) = 1$ and $H_1(x,s,q) = x$.

The first values are
$1,\ x,\ x^2 - (1-q)s,\ x^3 - (1-q)(q+2)sx,\ x^4 - (1-q)(q^2+2q+3)sx^2 + (1-q)^2(1+q+q^2)s^2, \cdots$.

Note that (5.2) implies that the coefficients of $(1-q)^k s^k$ are polynomials in $q$ with integer coefficients.

The polynomials $H_n(x,s,q)$ satisfy

$$H_n(x + (1-q)sD_q, s, q)1 = x^n. \tag{5.4}$$

This follows by induction since

$$H_{n+1}(x+(1-q)sD_q, s, q)1 = \left(x+(1-q)sD_q\right)H_n(x+(1-q)sD_q, s, q)1 + \left(q^n - 1\right)sH_{n-1}(x+(1-q)sD_q, s, q)1$$
$$= \left(x+(1-q)sD_q\right)x^n + \left(q^n - 1\right)sx^{n-1} = x^{n+1} + \left(1-q^n\right)sx^{n-1} + \left(q^n - 1\right)sx^{n-1} = x^{n+1}.$$

Therefore the umbral inverse sequence is given by

$$h_n(x,s,q) = \left(x + (1-q)sD_q\right)^n 1. \tag{5.5}$$

## 5.2. Relationship with Rogers-Szegö polynomials

The $q$ – Hermite polynomials $H_n(x,s,q)$ are intimately related with the Rogers-Szegö polynomials. Comparing (2.9) with (5.3) we see that

$$H_n(x,s,q) = R_n(\alpha,\beta,q) = \sum_{k=0}^{n}\begin{bmatrix}n\\k\end{bmatrix}\alpha^k\beta^{n-k} \tag{5.6}$$

with $\alpha = \dfrac{x+\sqrt{x^2-4s}}{2}$ und $\beta = \dfrac{x-\sqrt{x^2-4s}}{2}$.

In general $H_n(x,s,q)$ has no simple evaluations, but if we choose $(x,s,q) = \left(\sqrt{q}+\dfrac{1}{\sqrt{q}},1,q^2\right)$ we get $\alpha = \sqrt{q}$ and $\beta = \dfrac{1}{\sqrt{q}}$.

By (2.11) this gives

$$H_n\left(\sqrt{q}+\frac{1}{\sqrt{q}},1,q^2\right) = \frac{1}{\sqrt{q^n}}\sum_{k=0}^{n}\begin{bmatrix}n\\k\end{bmatrix}_{q^2} q^k = \frac{(1+q)(1+q^2)\cdots(1+q^n)}{\sqrt{q^n}}. \tag{5.7}$$

For $(x,s,q) = (1-q,-q,q)$ we get

$$H_n(1-q,-q,q) = \sum_{k=0}^{n}(-q)^k\begin{bmatrix}n\\k\end{bmatrix} = (1-q)(1-q^3)\cdots\left(1-q^{2\left\lfloor\frac{n+1}{2}\right\rfloor-1}\right). \tag{5.8}$$

From (2.16) we see that

$$\sum_{n\geq 0}\frac{H_n(x,s,q)}{[n]!}t^n = e\left(\frac{x+\sqrt{x^2-4s}}{2}t\right)e\left(\frac{x-\sqrt{x^2-4s}}{2}t\right). \tag{5.9}$$

For $q=0$ this reduces to the generating function of the Fibonacci polynomials

$$\sum_{n\geq 0}F_{n+1}(x,-s)t^n = \frac{1}{\left(1-\dfrac{x+\sqrt{x^2-4s}}{2}t\right)\left(1-\dfrac{x-\sqrt{x^2-4s}}{2}t\right)} = \frac{1}{1-xt+st^2}.$$

It is easy to verify that (5.9) satisfies the differential equation (5.1).

Note that $\dfrac{\partial}{\partial t}f(t)g(t) = g(t)\dfrac{\partial}{\partial t}f(t) + f(qt)\dfrac{\partial}{\partial t}g(t)$ and that $\dfrac{\partial}{\partial t}f(at) = af(t)$ and $e(qt) = (1+(q-1)t)e(t)$.

Therefore

$$\frac{\partial}{\partial t}\sum_{n\geq 0}\frac{H_n(x,s,q)}{[n]!}t^n = e\left(\frac{x-\sqrt{x^2-4s}}{2}t\right)\frac{\partial}{\partial t}e\left(\frac{x+\sqrt{x^2-4s}}{2}t\right) + e\left(q\frac{x+\sqrt{x^2-4s}}{2}t\right)\frac{\partial}{\partial t}e\left(\frac{x-\sqrt{x^2-4s}}{2}t\right)$$

$$=\left(\frac{x+\sqrt{x^2-4s}}{2}+\frac{x-\sqrt{x^2-4s}}{2}\left(1-(1-q)\frac{x+\sqrt{x^2-4s}}{2}t\right)\right)e\left(\frac{x+\sqrt{x^2-4s}}{2}t\right)e\left(\frac{x-\sqrt{x^2-4s}}{2}t\right)$$

$$=(x-(1-q)st)\sum_{n\geq 0}\frac{H_n(x,s,q)}{[n]!}t^n.$$

By (2.25) we can write (5.9) as

$$\sum_{n\geq 0}\frac{H_n(x,s,q)}{(q;q)_n}t^n = \frac{1}{\left(\frac{x+\sqrt{x^2-4s}}{2}t;q\right)_\infty \left(\frac{x-\sqrt{x^2-4s}}{2}t;q\right)_\infty} = \frac{1}{\prod_{n=0}^{\infty}\left(1-q^n xt + q^{2n}st^2\right)}. \quad (5.10)$$

**5.3. Connection with Fibonacci and Lucas polynomials**

The $q$ – Hermite polynomials are also related to the Fibonacci and Lucas polynomials.

To this end let us take a closer look at formula (5.6).

For odd $n$ we have

$$\sum_{k=0}^{n}\begin{bmatrix}n\\k\end{bmatrix}\alpha^k\beta^{n-k} = \sum_{k=0}^{\left\lfloor\frac{n}{2}\right\rfloor}\begin{bmatrix}n\\k\end{bmatrix}\left(\alpha^k\beta^{n-k}+\alpha^{n-k}\beta^k\right) = \sum_{k=0}^{\left\lfloor\frac{n}{2}\right\rfloor}\begin{bmatrix}n\\k\end{bmatrix}\alpha^k\beta^k\left(\beta^{n-2k}+\alpha^{n-2k}\right) = \sum_{k=0}^{\left\lfloor\frac{n}{2}\right\rfloor}\begin{bmatrix}n\\k\end{bmatrix}s^k\left(\beta^{n-2k}+\alpha^{n-2k}\right).$$

For $n=2m$ we get

$$\sum_{k=0}^{2m}\begin{bmatrix}2m\\k\end{bmatrix}\alpha^k\beta^{2m-k} = \begin{bmatrix}2m\\m\end{bmatrix}\alpha^m\beta^m + \sum_{k=0}^{m-1}\begin{bmatrix}2m\\k\end{bmatrix}\left(\alpha^k\beta^{2m-k}+\alpha^{2m-k}\beta^k\right) = \sum_{k=0}^{m-1}\begin{bmatrix}2m\\k\end{bmatrix}s^k\left(\beta^{2m-2k}+\alpha^{2m-2k}\right) + \begin{bmatrix}2m\\m\end{bmatrix}s^m$$

Therefore we get

$$\sum_{k=0}^{\left\lfloor\frac{n}{2}\right\rfloor}\begin{bmatrix}n\\k\end{bmatrix}s^k L^*_{n-2k}(x,-s) = H_n(x,s,q). \quad (5.11)$$

For $q=1$ this reduces to (3.5).

For $q=0$ this reduces to $\sum_{k=0}^{\left\lfloor\frac{n}{2}\right\rfloor}s^k L^*_{n-2k}(x,-s) = F_{n+1}(x,-s).$

The latter fact is rather trivial because $L_n(x,-s) = F_{n+1}(x,-s) - sF_{n-1}(x,-s)$ reduces the left-hand side to a telescoping sum.

Expressing the Lucas polynomials by Fibonacci polynomials we get

$$H_n(x,s,q) = \sum_{k=0}^{\lfloor n/2 \rfloor} \begin{bmatrix} n \\ k \end{bmatrix} s^k L^*_{n-2k}(x,-s) = \sum_{k=0}^{\lfloor n/2 \rfloor} \begin{bmatrix} n \\ k \end{bmatrix} s^k F_{n+1-2k}(x,-s) - \sum_{k=0}^{\lfloor n/2 \rfloor - 1} \begin{bmatrix} n \\ k \end{bmatrix} s^{k+1} F_{n-1-2k}(x,-s)$$

$$= \sum_{k=0}^{\lfloor n/2 \rfloor} \begin{bmatrix} n \\ k \end{bmatrix} s^k F_{n+1-2k}(x,-s) - \sum_{k=0}^{\lfloor n/2 \rfloor} \begin{bmatrix} n \\ k-1 \end{bmatrix} s^k F_{n+1-2k}(x,-s) = \sum_{k=0}^{\lfloor n/2 \rfloor} \left( \begin{bmatrix} n \\ k \end{bmatrix} - \begin{bmatrix} n \\ k-1 \end{bmatrix} \right) s^k F_{n+1-2k}(x,-s).$$

Thus we have

$$\sum_{k=0}^{\lfloor n/2 \rfloor} \left( \begin{bmatrix} n \\ k \end{bmatrix} - \begin{bmatrix} n \\ k-1 \end{bmatrix} \right) s^k F_{n+1-2k}(x,-s) = H_n(x,s,q) \tag{5.12}$$

This gives an explicit evaluation of the coefficients of $H_n(x,s,q)$:

$$H_n(x,s,q) = \sum_{k=0}^{\lfloor n/2 \rfloor} h(n,k,q) s^k x^{n-2k} \tag{5.13}$$

with

$$h(n,k,q) = \sum_{j=0}^{k} \left( \begin{bmatrix} n \\ j \end{bmatrix} - \begin{bmatrix} n \\ j-1 \end{bmatrix} \right) \binom{n-k-j}{k-j} (-1)^{k-j}. \tag{5.14}$$

We know that $\dfrac{h(n,k,q)}{(1-q)^k}$ is a polynomial in $q$ with integer coefficients. For $k=1$ we get $h(n,1,q) = -(1-q)\left(q^{n-2} + 2q^{n-3} + 3q^{n-4} + \cdots + (n-1)\right)$ but for $k>1$ no obvious simple expression occurs.

### 5.4. Associated inverse relations

If in (5.11) we replace $x$ by $x + (1-q)sD_q$ we get by (5.4)

$$\sum_{k=0}^{\lfloor n/2 \rfloor} \begin{bmatrix} n \\ k \end{bmatrix} s^k L^*_{n-2k}(x + (1-q)sD_q, -s)1 = x^n. \tag{5.15}$$

By replacing $x$ by $x + (1-q)sD_q$ in (5.12) we get the analogous formula

$$\sum_{k=0}^{\lfloor n/2 \rfloor} \left( \begin{bmatrix} n \\ k \end{bmatrix} - \begin{bmatrix} n \\ k-1 \end{bmatrix} \right) s^k F_{n+1-2k}(x + (1-q)sD_q, -s) = x^n. \tag{5.16}$$

It comes as a big surprise that both $F_{n+1}(x+(1-q)sD_q,-s)1$ and $L_n^*(x+(1-q)sD_q,-s)1$ have very simple expressions (cf. [7]):

$$F_{n+1}(x+(1-q)sD_q,-s)1 = F_{n+1}(x,-s,q) = \sum_{k=0}^{\lfloor \frac{n}{2} \rfloor} \begin{bmatrix} n-k \\ k \end{bmatrix} q^{\binom{k+1}{2}}(-s)^k x^{n-2k} \quad (5.17)$$

and for $n > 0$

$$L_n^*(x+(1-q)sD_q,-s)1 = L_n(x,-s,q) = \sum_{k=0}^{\lfloor \frac{n}{2} \rfloor} q^{\binom{k}{2}} \frac{[n]}{[n-k]} \begin{bmatrix} n-k \\ k \end{bmatrix}(-s)^k x^{n-2k}. \quad (5.18)$$

We consider first the case of the Fibonacci polynomials.

Let $F_n(x,s,q) = F_n(x+(q-1)sD_q,s)1 = \big(x+(q-1)sD_q\big)F_{n-1}(x,s,q) + sF_{n-2}(x,s,q)$

The first terms are

$0, 1, x, x^2 + qs, x^3 + q[2]sx, x^4 + q[3]sx^2 + q^3s^2, x^5 + q[4]sx^3 + q^3[3]s^2x,\cdots$.

This leads us to guess that

$$F_n(x,s,q) = \begin{bmatrix} n-1 \\ 0 \end{bmatrix} x^{n-1} + q\begin{bmatrix} n-2 \\ 1 \end{bmatrix} sx^{n-3} + q^3\begin{bmatrix} n-3 \\ 2 \end{bmatrix} s^2x^{n-5} + \cdots = \sum_{k=0}^{\lfloor \frac{n-1}{2} \rfloor} \begin{bmatrix} n-1-k \\ k \end{bmatrix} q^{\binom{k+1}{2}} s^k x^{n-1-2k}.$$

Since $F_n(x,s,q) = \big(x+(q-1)sD_q\big)F_{n-1}(x,s,q) + sF_{n-2}(x,s,q)$ by comparing coefficients this conjecture is equivalent to

$$\begin{bmatrix} n-1-k \\ k \end{bmatrix} q^{\binom{k+1}{2}} = \begin{bmatrix} n-2-k \\ k \end{bmatrix} q^{\binom{k+1}{2}} + \begin{bmatrix} n-1-k \\ k-1 \end{bmatrix} q^{\binom{k}{2}}(q^{n-2k}-1) + q^{\binom{k}{2}}\begin{bmatrix} n-2-k \\ k-1 \end{bmatrix}$$

$$= q^{\binom{k}{2}}\left(\begin{bmatrix} n-1-k \\ k \end{bmatrix} + q^{n-2k}\begin{bmatrix} n-1-k \\ k-1 \end{bmatrix}\right) - \begin{bmatrix} n-1-k \\ k-1 \end{bmatrix} q^{\binom{k}{2}} = q^{\binom{k}{2}}\left(\begin{bmatrix} n-k \\ k \end{bmatrix} - \begin{bmatrix} n-1-k \\ k-1 \end{bmatrix}\right)$$

which is obviously true.

Now we get

$$L_n(x,s,q) = L_n\big(x+(q-1)sD_q,s\big)1 = F_{n+1}\big(x+(q-1)sD_q,s\big)1 + sF_{n-1}\big(x+(q-1)sD_q,s\big)1$$

$$= F_{n+1}(x,s,q) + sF_{n-1}(x,s,q) = \sum_{k=0}^{\lfloor \frac{n}{2} \rfloor} \begin{bmatrix} n-k \\ k \end{bmatrix} q^{\binom{k+1}{2}} s^k x^{n-2k} + \sum_{k=0}^{\lfloor \frac{n-2}{2} \rfloor} \begin{bmatrix} n-2-k \\ k \end{bmatrix} q^{\binom{k+1}{2}} s^{k+1} x^{n-2-2k}$$

$$= \sum_{k=0}^{\lfloor \frac{n}{2} \rfloor}\left(\begin{bmatrix} n-k \\ k \end{bmatrix} q^k + \begin{bmatrix} n-1-k \\ k-1 \end{bmatrix}\right) q^{\binom{k}{2}} s^k x^{n-2k} = \sum_{k=0}^{\lfloor \frac{n}{2} \rfloor} q^{\binom{k}{2}} \frac{[n]}{[n-k]} \begin{bmatrix} n-k \\ k \end{bmatrix} s^k x^{n-2k}.$$

Thus we have

$$L_n(x,s,q) = \sum_{k=0}^{\lfloor \frac{n}{2} \rfloor} q^{\binom{k}{2}} \frac{[n]}{[n-k]} \begin{bmatrix} n-k \\ k \end{bmatrix} s^k x^{n-2k} \tag{5.19}$$

and

$$\sum_{k=0}^{\lfloor \frac{n}{2} \rfloor} \begin{bmatrix} n \\ k \end{bmatrix} s^k L^*_{n-2k}(x,-s,q) = x^n. \tag{5.20}$$

This implies the inversion formula ([4], Theorem 6):

Let

$$\sum_{k=0}^{\lfloor \frac{n}{2} \rfloor} \begin{bmatrix} n \\ k \end{bmatrix} s^k a(n-2k) = b(n). \tag{5.21}$$

Then

$$a(n) = \sum_{k=0}^{\lfloor \frac{n}{2} \rfloor} q^{\binom{k}{2}} \frac{[n]}{[n-k]} \begin{bmatrix} n-k \\ k \end{bmatrix} (-s)^k b(n-2k) \tag{5.22}$$

and vice versa.

In the same way we obtain from the identities (5.17) and

$$\sum_{k=0}^{\lfloor \frac{n}{2} \rfloor} \left( \begin{bmatrix} n \\ k \end{bmatrix} - \begin{bmatrix} n \\ k-1 \end{bmatrix} \right) s^k F_{n+1-2k}(x,-s,q) = x^n \tag{5.23}$$

the equivalence of

$$\sum_{k=0}^{\lfloor \frac{n}{2} \rfloor} \left( \begin{bmatrix} n \\ k \end{bmatrix} - \begin{bmatrix} n \\ k-1 \end{bmatrix} \right) s^k a(n-2k) = b(n) \tag{5.24}$$

and

$$a(n) = \sum_{k=0}^{\lfloor \frac{n}{2} \rfloor} \begin{bmatrix} n-k \\ k \end{bmatrix} q^{\binom{k+1}{2}} (-s)^k b(n-2k). \tag{5.25}$$

This is Theorem 7 of [4].

If we define $\Lambda_{q,F}$ by $\Lambda_{q,F}\left(F_{n+1}(x,-1,q)\right) = [n=0]$ we get

$$\Lambda_{q,F}\left(x^{2n}\right) = q^n \frac{1}{[n+1]} \begin{bmatrix} 2n \\ n \end{bmatrix}. \tag{5.26}$$

If we define $\Lambda_{q,L}$ by $\Lambda_{q,L}\left(L_n(x,-1,q)\right) = [n=0]$ we get

$$\Lambda_{q,L}\left(x^{2n}\right) = \begin{bmatrix} 2n \\ n \end{bmatrix}. \tag{5.27}$$

Note that neither the $q$-Fibonacci polynomials nor the $q$-Lucas polynomials are orthogonal with respect to these linear functionals.

If we choose $a(n) = L_n^*(x,-s)$ in (5.21) then (5.11) implies

$$L_n^*(x,-s) = \sum_{k=0}^{\lfloor n/2 \rfloor} q^{\binom{k}{2}} \frac{[n]}{[n-k]} \begin{bmatrix} n-k \\ k \end{bmatrix} (-s)^k H_{n-2k}(x,s) \tag{5.28}$$

and (5.12) implies for $a(n) = F_{n+1}(x,-s)$

$$F_{n+1}(x,-s) = \sum_{k=0}^{\lfloor n/2 \rfloor} \begin{bmatrix} n-k \\ k \end{bmatrix} q^{\binom{k+1}{2}} (-s)^k H_{n-2k}(x,s,q). \tag{5.29}$$

## 5.5. Moments of the q-Hermite polynomials

In [10] we introduced the polynomials $h_n(x,s,q) = \left(x + (1-q)sD_q\right)^n 1$.
If in (3.4) and (3.5) we replace $x$ by $x + (1-q)sD_q$ we get

$$\sum_{k=0}^{\lfloor n/2 \rfloor} \binom{n}{k} s^k L_{n-2k}^*(x+(1-q)sD_q,-s)1 = \left(x+(1-q)sD_q\right)^n 1 \text{ or}$$

$$h_n(x,s,q) = \sum_{k=0}^{\lfloor n/2 \rfloor} \binom{n}{k} s^k L_{n-2k}^*(x,-s,q) = \sum_{k=0}^{\lfloor n/2 \rfloor} \left(\binom{n}{k} - \binom{n}{k-1}\right) s^k F_{n+1-2k}(x,-s,q). \tag{5.30}$$

From this we see that $h_n(x,s,q) = \sum_{k=0}^{\lfloor n/2 \rfloor} c(n,k,q) s^k x^{n-2k}$

with

$$c(n,k,q) = \sum_{j=0}^{k} \left(\binom{n}{j} - \binom{n}{j-1}\right)(-1)^{k-j} q^{\binom{k-j+1}{2}} \begin{bmatrix} n-j-k \\ k-j \end{bmatrix} = \sum_{j=0}^{k} (-1)^{k-j} \binom{n}{j} q^{\binom{k-j}{2}} \frac{[n-2j]}{[n-j-k]} \begin{bmatrix} n-j-k \\ k-j \end{bmatrix}.$$

Again $\dfrac{c(n,k,q)}{(1-q)^k}$ is a polynomial in $q$ with integer coefficients.

By (5.4) we have $H_n(x+(1-q)sD_q,s,q)1 = x^n$.

Therefore we get the **Touchard-Riordan formula** (cf. e.g. [10],[15])

$$\Lambda_{H(x,s,q)}\left(x^{2n}\right) = c(2n,n,q) = s^n \sum_{j=0}^{n} \left(\binom{2n}{n-j} - \binom{2n}{n-j-1}\right)(-1)^j q^{\binom{j+1}{2}}. \tag{5.31}$$

# 6. The continuous q-Hermite polynomials $H_n(x|q)$.

## 6.1. Some simple properties

The continuous $q$ – Hermite polynomials $H_n(x|q) = H_n(2x,1,q)$ satisfy the recurrence

$$H_{n+1}(x|q) = 2xH_n(x|q) - (1-q^n)H_{n-1}(x|q) \tag{6.1}$$

with initial values $H_0(x|q) = 1$ and $H_1(x|q) = 2x$.
The first values are

$$1,\ 2x,\ 4x^2 + q - 1,\ 8x^3 + 2(q^2 + q - 2)x,\ 16x^4 + 4(q^3 + q^2 + q - 3)x^2 + q^4 - q^3 - q + 1,\ \cdots.$$

The rescaled polynomials $\tilde{H}_n(x|q) = \dfrac{H_n(x\sqrt{1-q}\,|q^2)}{\sqrt{(1-q)^n}}$ are a $q$ – analogue of the physicists'
Hermite polynomials since they satisfy $\tilde{H}_{n+1}(x|q) = 2x\tilde{H}_n(x|q) - [2n]_q \tilde{H}_{n-1}(x|q)$ and therefore
also $\lim_{q \to 1} \tilde{H}_n(x|q) = \mathbf{H}_n(x)$.

There are many connections with trigonometric functions.

For $H_n(x|q)$ we have $\alpha = x + \sqrt{x^2 - 1}$ and $\beta = x - \sqrt{x^2 - 1}$. If we set $x = \cos\vartheta$ then (5.6) gives the well-known formula

$$H_n(x|q) = \sum_{k=0}^{n} \begin{bmatrix} n \\ k \end{bmatrix} e^{i(n-2k)\vartheta}. \tag{6.2}$$

The generating function reduces to

$$\sum_{n \geq 0} \frac{H_n(x|q)}{(q;q)_n} t^n = \frac{1}{\left((x+\sqrt{x^2-1})t;q\right)_\infty \left((x-\sqrt{x^2-1})t;q\right)_\infty} = \frac{1}{\prod_{n=0}^{\infty}(1 - 2q^n xt + q^{2n}t^2)}$$

or equivalently to

$$\sum_{n \geq 0} \frac{H_n(\cos\vartheta|q)}{(q;q)_n} t^n = \frac{1}{(e^{i\vartheta}t;q)_\infty (e^{-i\vartheta}t;q)_\infty} = \frac{1}{\left|(e^{i\vartheta}t;q)_\infty\right|^2}. \tag{6.3}$$

For the continuous $q$ – Hermite polynomials the role of Fibonacci and Lucas polynomials is played by Chebyshev polynomials.

Here (5.12) takes the form

$$\sum_{k=0}^{\lfloor \frac{n}{2} \rfloor} \left( \begin{bmatrix} n \\ k \end{bmatrix} - \begin{bmatrix} n \\ k-1 \end{bmatrix} \right) U_{n-2k}(x) = H_n(x|q) \tag{6.4}$$

and the inverse

$$U_n(x) = \sum_{k=0}^{\lfloor \frac{n}{2} \rfloor} \begin{bmatrix} n-k \\ k \end{bmatrix} q^{\binom{k+1}{2}} (-1)^k H_{n-2k}(x \mid q). \tag{6.5}$$

On the other hand we get

$$H_n(x \mid q) = \sum_{k=0}^{n} \begin{bmatrix} n \\ k \end{bmatrix} T_{n-2k}(x), \tag{6.6}$$

if we write $T_{-n}(x) = T_n(x)$.

The inverse is

$$T_n(x) = \frac{1}{2} \sum_{k=0}^{\lfloor \frac{n}{2} \rfloor} q^{\binom{k}{2}} \frac{[n]}{[n-k]} \begin{bmatrix} n-k \\ k \end{bmatrix} (-1)^k H_{n-2k}(x \mid q). \tag{6.7}$$

### 6.2. Orthogonality and the associated probability measure

Let now $\Lambda_{H,q}$ be the linear functional defined by $\Lambda_{H,q}(H_n(x \mid q)) = [n = 0]$. Since the polynomials $H_n(x \mid q)$ satisfy the 3-term recurrence (6.1) they are orthogonal with respect to $\Lambda_{H,q}$.
More precisely we have

$$\Lambda_{H,q}(H_n(x \mid q) H_m(x \mid q)) = (q;q)_n [n = m]. \tag{6.8}$$

Since $H_n(x \mid q) = (2x)^n + \cdots$ it suffices to show that
$\Lambda_{H,q}((2x)^n H_n(x \mid q)) = (q;q)_n$.
This follows by multiplying (6.1) by $(2x)^{n-1}$ and applying $\Lambda_{H,q}$.

Now we want to determine the associated probability measure (cf. [1], [2], [11], [13], [18], [19]).

By (6.5) we see that $\Lambda_{H,q}(U_{2n+1}) = 0$ and

$$\Lambda_{H,q}(U_{2n}) = (-1)^n q^{\binom{n+1}{2}}. \tag{6.9}$$

From (6.7) we get $\Lambda_{H,q}(T_{2n+1}) = 0$ and

$$\Lambda_{H,q}(T_{2n}) = \frac{(-1)^n}{2} q^{\binom{n}{2}} (1 + q^n). \tag{6.10}$$

Formula (6.9) implies

$$\Lambda_{H,q}(f(x)) = \int_{-1}^{1} f(x) w(x;q) dx \tag{6.11}$$

where

$$w(x;q) = \frac{2}{\pi} \sqrt{1-x^2} \sum_{k=0}^{\infty} (-1)^k q^{\binom{k+1}{2}} U_{2k}(x). \tag{6.12}$$

For the orthogonality of the Chebyshev polynomials of the second kind gives

$$\int_{-1}^{1} U_n(x) U_k(x) \sqrt{1-x^2}\, dx = \frac{\pi}{2}[n=k] \tag{6.13}$$

and therefore

$$\int_{-1}^{1} U_{2n+1}(x) w(x;q)\, dx = \frac{2}{\pi} \sum_{k=0}^{\infty} (-1)^k q^{\binom{k+1}{2}} \int_{-1}^{1} U_{2n+1}(x) U_{2k}(x) \sqrt{1-x^2}\, dx = 0 = \Lambda_{H,q}(U_{2n+1})$$

and

$$\int_{-1}^{1} U_{2n}(x) w(x;q)\, dx = \frac{2}{\pi} \sum_{k=0}^{\infty} (-1)^k q^{\binom{k+1}{2}} \int_{-1}^{1} U_{2n}(x) U_{2k}(x) \sqrt{1-x^2}\, dx = (-1)^n q^{\binom{n+1}{2}} = \Lambda_{H,q}(U_{2n}).$$

Therefore (6.11) holds for each polynomial $f(x)$.
Formula (6.12) has first been obtained by Wm. R. Allaway [1] and later also by Pawel J. Szablowski [18],[19].

In most papers this functional is written in the following form:

$$\Lambda_{H,q}(f(\cos\vartheta)) = \frac{(q;q)_\infty}{2\pi} \int_0^\pi f(\cos\vartheta)\left(e^{2i\vartheta};q\right)_\infty \left(e^{-2i\vartheta};q\right)_\infty d\vartheta. \tag{6.14}$$

Since $\left(e^{2i\vartheta};q\right)_\infty \left(e^{-2i\vartheta};q\right)_\infty = \left|\left(e^{2i\vartheta};q\right)_\infty\right|^2 \geq 0$ this defines a probability measure.

To show (6.14) let $x = \cos\vartheta$, $0 \leq \vartheta < \pi$. (6.11) gives
$$\Lambda_{H,q}(f(\cos\vartheta)) = \int_0^\pi f(\cos\vartheta) w(\cos\vartheta;q) \sin\vartheta\, d\vartheta$$

$$= \frac{2}{\pi} \int_0^\pi f(\cos\vartheta) \sin\vartheta \sum_{k=0}^{\infty} (-1)^k q^{\binom{k+1}{2}} \frac{\sin(2k+1)\vartheta}{\sin\vartheta} \sin\vartheta\, d\vartheta$$

$$= \frac{2}{\pi} \int_0^\pi f(\cos\vartheta) \sum_{k=0}^{\infty} (-1)^k q^{\binom{k+1}{2}} \sin(2k+1)\vartheta \sin\vartheta\, d\vartheta.$$

Now we apply Jacobi's triple product identity

$$\sum_{k\in\mathbb{Z}} (-1)^k q^{\binom{k}{2}} x^k = (x;q)_\infty \left(\frac{q}{x};q\right)_\infty (q;q)_\infty. \tag{6.15}$$

First write

$$2i \sum_{k=0}^{\infty} (-1)^k q^{\binom{k+1}{2}} \sin(2k+1)\vartheta = \sum_{k\in\mathbb{Z}} (-1)^k q^{\binom{k+1}{2}} e^{(2k+1)i\vartheta} = -e^{-i\vartheta} \sum_{k\in\mathbb{Z}} (-1)^k q^{\binom{k}{2}} \left(e^{2i\vartheta}\right)^k.$$

This implies

$$\sin\vartheta\sum_{k=0}^{\infty}(-1)^k q^{\binom{k+1}{2}}\sin(2k+1)\vartheta = -\frac{e^{-i\vartheta}\sin\vartheta}{2i}(e^{2i\vartheta};q)_\infty(qe^{-2i\vartheta};q)_\infty(q;q)_\infty$$

$$= \frac{1}{4}(e^{2i\vartheta};q)_\infty(e^{-2i\vartheta};q)_\infty(q;q)_\infty$$

and therefore (6.14).

We will now give another proof of (6.14) without using Jacobi's triple product identity.

Since $T_n(\cos\vartheta) = \cos n\vartheta$ identity (6.14) is equivalent with

$$\frac{(q;q)_\infty}{2\pi}\int_0^\pi \cos n\vartheta\, (e^{2i\vartheta};q)_\infty(e^{-2i\vartheta};q)_\infty\, d\vartheta = \Lambda_{H,q}(T_n) \tag{6.16}$$

for all $n \in \mathbb{N}$.

Since $(e^{2i\vartheta};q)_\infty(e^{-2i\vartheta};q)_\infty = |(e^{2i\vartheta};q)_\infty|^2$ is real we have

$$I_k := \int_0^\pi \cos k\vartheta\,(e^{2i\vartheta};q)_\infty(e^{-2i\vartheta};q)_\infty\, d\vartheta = \mathrm{Re}\int_0^\pi e^{ik\vartheta}(e^{2i\vartheta};q)_\infty(e^{-2i\vartheta};q)_\infty\, d\vartheta.$$

Now observe that

$$(e^{2i\vartheta};q)_\infty(e^{-2i\vartheta};q)_\infty = \sum_j \frac{(-1)^j q^{\binom{j}{2}} e^{2ij\vartheta}}{(q;q)_j}\sum_\ell \frac{(-1)^\ell q^{\binom{\ell}{2}} e^{-2i\ell\vartheta}}{(q;q)_\ell} = \sum_{j,\ell}(-1)^{j+\ell}\frac{q^{\binom{j}{2}+\binom{\ell}{2}}}{(q;q)_j(q;q)_\ell} e^{2i(j-\ell)\vartheta}$$

and that $\int_0^\pi e^{i(2k+1)\vartheta} d\vartheta = \dfrac{2i}{2k+1}$ is purely imaginary. Therefore $I_k = 0$ if $k$ is odd.

On the other hand we have

$$I_{2k} := \int_0^\pi e^{2ik\vartheta}(e^{2i\vartheta};q)_\infty(e^{-2i\vartheta};q)_\infty\, d\vartheta = \sum_{j,\ell}(-1)^{j+\ell}\frac{q^{\binom{j}{2}+\binom{\ell}{2}}}{(q;q)_j(q;q)_\ell}\int_0^\pi e^{i(2j-2\ell+2k)\vartheta} d\vartheta$$

$$= \pi\sum_j \frac{(-1)^k q^{\binom{j}{2}+\binom{j+k}{2}}}{(q;q)_j(q;q)_{j+k}} = \pi\frac{(-1)^k q^{\binom{k}{2}}}{(q;q)_k}\sum_j \frac{q^{j^2-j}}{(q;q)_j(q^{k+1};q)_j}(q^k)^j$$

Therefore because of (6.10) it suffices to prove that

$$\sum_j \frac{q^{j(j-1)}}{(q;q)_j(qq^k;q)_j}(q^k)^j = \frac{(1+q^k)}{(q^{k+1};q)_\infty}. \tag{6.17}$$

We will prove more generally

$$\sum_{j=0}^{\infty} \frac{q^{j^2-j} z^j \left(q^{j+1} z; q\right)_\infty}{(q;q)_j} = 1 + z, \qquad (6.18)$$

which for $z = q^k$ coincides with (6.17).

The left-hand side is

$$\sum_{j=0}^{\infty} \frac{q^{j^2-j} z^j}{(q;q)_j} \sum_{k=0}^{\infty} (-1)^k q^{\binom{k}{2}} \frac{\left(q^{j+1} z\right)^k}{(q;q)_k} = \sum_{n \geq 0} \frac{z^n}{(q;q)_n} \sum_{k=0}^{n} (-1)^k \begin{bmatrix} n \\ k \end{bmatrix} q^{\binom{k}{2}+k(n-k+1)+(n-k)^2-(n-k)}$$

$$= \sum_{n \geq 0} \frac{q^{n^2-n}}{(q;q)_n} z^n \sum_{k=0}^{n} (-1)^k \begin{bmatrix} n \\ k \end{bmatrix} q^{\frac{k^2+3k}{2}-nk}.$$

Let us calculate $w(n) = \sum_{k=0}^{n} (-1)^k \begin{bmatrix} n \\ k \end{bmatrix} q^{\frac{k^2+3k}{2}-nk}$.

We have $w(0) = 1$ and $w(1) = 1 - q^{2-1} = 1 - q$.

For $n \geq 2$

$$w(n) = \sum_{k=0}^{n} (-1)^k \begin{bmatrix} n \\ k \end{bmatrix} q^{\frac{k^2-k}{2}} \left(q^{-(n-2)}\right)^k = \left(1 - q^{-(n-2)}\right)\left(1 - q^{-(n-3)}\right) \cdots (1-1)(1-q) = 0.$$

We have thus proved that

$$\frac{(q;q)_\infty}{2\pi} \int_0^\pi \cos 2n\vartheta \left(e^{2i\vartheta}; q\right)_\infty \left(e^{-2i\vartheta}; q\right)_\infty d\vartheta = \frac{(-1)^n}{2} q^{\binom{n}{2}} (1+q^n) = \Lambda_{H,q}(T_{2n})$$

and therefore (6.16).

### 6.3. The Askey-Wilson q-differentiation operator

The Hermite polynomials $\mathbf{H}_n(x)$ satisfy $D\mathbf{H}_n(x) = 2n\mathbf{H}_{n-1}(x)$. Askey and Wilson [3] have found an analogue of this relation for the continuous $q$ – Hermite polynomials.

Define $\delta_q F(x)$ for a function $F(x) = f\left(e^{i\vartheta}\right)$ of $x = \cos \vartheta$ by

$$\delta_q F(x) = \delta_q f\left(e^{i\vartheta}\right) = f\left(q^{\frac{1}{2}} e^{i\vartheta}\right) - f\left(q^{-\frac{1}{2}} e^{i\vartheta}\right).$$

For $F(x) = \cos \vartheta = \dfrac{e^{i\vartheta} + e^{-i\vartheta}}{2}$ this gives

$$\delta_q x = \delta_q \frac{e^{i\vartheta} + e^{-i\vartheta}}{2} = \frac{1}{2}\left(\sqrt{q} e^{i\vartheta} + \frac{1}{\sqrt{q}} e^{i\vartheta} - \frac{1}{\sqrt{q}} e^{-i\vartheta} - \sqrt{q} e^{-i\vartheta}\right) = i\left(\sqrt{q} - \frac{1}{\sqrt{q}}\right) \sin \vartheta.$$

Then the **Askey-Wilson operator** $\Delta_q$ is defined by

$$\Delta_q F(x) = \Delta_q f\left(e^{i\vartheta}\right) = \frac{\delta_q f\left(e^{i\vartheta}\right)}{\delta_q \cos\vartheta} = \frac{\delta_q f\left(e^{i\vartheta}\right)}{i\left(q^{\frac{1}{2}} - q^{-\frac{1}{2}}\right)\sin\vartheta}. \tag{6.19}$$

To get a feeling for this operator let us remark that

$$\Delta_q T_n(x) = \frac{[n]}{\sqrt{q^{n-1}}} U_{n-1}(x) \tag{6.20}$$

for $x = \cos\vartheta$.

This follows immediately from $T_n(\cos\vartheta) = \cos n\vartheta = \frac{1}{2}\left(e^{in\vartheta} + e^{-in\vartheta}\right)$ since

$$\delta_q \frac{1}{2}\left(e^{in\vartheta} + e^{-in\vartheta}\right) = \frac{1}{2}\left(q^{\frac{n}{2}}e^{in\vartheta} - q^{-\frac{n}{2}}e^{in\vartheta} + q^{-\frac{n}{2}}e^{-in\vartheta} + q^{\frac{n}{2}}e^{-in\vartheta}\right) = i\left(q^{\frac{n}{2}} - q^{-\frac{n}{2}}\right)\sin n\vartheta$$

and therefore

$$\Delta_q \frac{1}{2}\left(e^{in\vartheta} + e^{-in\vartheta}\right) = \frac{i\left(q^{\frac{n}{2}} - q^{-\frac{n}{2}}\right)\sin n\vartheta}{i\left(q^{\frac{1}{2}} - q^{-\frac{1}{2}}\right)\sin\vartheta} = \frac{[n]}{q^{\frac{n-1}{2}}}\frac{\sin n\vartheta}{\sin\vartheta} = \frac{[n]}{q^{\frac{n-1}{2}}} U_{n-1}(\cos\vartheta).$$

The Askey-Wilson operator applied to the continuous $q$ – Hermite polynomials gives

$$\Delta_q H_n(x \mid q) = 2q^{-\frac{n-1}{2}}[n]_q H_{n-1}(x \mid q) \tag{6.21}$$

for $x = \cos\vartheta$.

This is equivalent with

$$\frac{q^{\frac{n}{2}}}{1-q^n}\Delta_q H_n(\cos\vartheta \mid q) = \frac{2q^{\frac{1}{2}}}{1-q} H_{n-1}(\cos\vartheta \mid q)$$

or

$$\frac{q^{\frac{n}{2}}}{1-q^n}\delta_q H_n(\cos\vartheta \mid q) = -\left(e^{i\vartheta} - e^{-i\vartheta}\right) H_{n-1}(\cos\vartheta \mid q).$$

The last identity follows from

$$\frac{q^{\frac{n}{2}}}{1-q^n}\sum_k \begin{bmatrix} n \\ k \end{bmatrix}\left(q^{\frac{n-2k}{2}} - q^{-\frac{n-2k}{2}}\right)e^{i\vartheta(n-2k)} = \sum_k \begin{bmatrix} n \\ k \end{bmatrix}\frac{q^{n-k}-q^k}{1-q^n}e^{i\vartheta(n-2k)}$$

$$= \sum_k \begin{bmatrix} n \\ k \end{bmatrix}\frac{1-q^k}{1-q^n}e^{i\vartheta(n-2k)} + \sum_k \begin{bmatrix} n \\ k \end{bmatrix}\frac{q^{n-k}-1}{1-q^n}e^{i\vartheta(n-2k)} = \sum_k \begin{bmatrix} n-1 \\ k-1 \end{bmatrix}e^{i\vartheta(n-2k)} - \sum_k \begin{bmatrix} n-1 \\ k \end{bmatrix}e^{i\vartheta(n-2k)}$$

$$= e^{-i\vartheta}\sum_k \begin{bmatrix} n-1 \\ k \end{bmatrix}e^{i\vartheta(n-1-2k)} - e^{i\vartheta}\sum_k \begin{bmatrix} n-1 \\ k \end{bmatrix}e^{i\vartheta(n-1-2k)} = -\left(e^{i\vartheta} - e^{-i\vartheta}\right)H_{n-1}(\cos\vartheta \mid q).$$

A $q$-analogue of (1.18) is

$$H_n(x \mid q) = \left(2x - \frac{1-q}{2}q^{\frac{n-2}{2}}\Delta_q\right)\left(2x - \frac{1-q}{2}q^{\frac{n-3}{2}}\Delta_q\right)\cdots\left(2x - \frac{1-q}{2}\Delta_q\right)1. \tag{6.22}$$

This follows from (6.1) and (6.21).

The polynomials $\mathbf{H}_n(x)$ satisfy $De^{-x^2}\mathbf{H}_n(x) = -2xe^{-x^2}\mathbf{H}_n(x) + e^{-x^2}2n\mathbf{H}_{n-1}(x) = -e^{-x^2}\mathbf{H}_{n+1}(x)$.

An analogous formula holds for the continuous $q$-Hermite polynomials, where the role of $e^{-x^2}$ is played by the function $v(\vartheta;q) = \dfrac{\left(e^{2i\vartheta};q\right)_\infty \left(e^{-2i\vartheta};q\right)_\infty}{\sin\vartheta}$.

$$\Delta_q \frac{\left(e^{2i\vartheta};q\right)_\infty \left(e^{-2i\vartheta};q\right)_\infty}{\sin\vartheta} H_n(\cos\vartheta \mid q) = \frac{-2}{1-q}\frac{\left(e^{2i\vartheta};q\right)_\infty \left(e^{-2i\vartheta};q\right)_\infty}{\sin\vartheta}\frac{H_{n+1}(\cos\vartheta \mid q)}{\sqrt{q^n}}. \tag{6.23}$$

It suffices to show that

$$\delta_q \frac{\left(e^{2i\vartheta};q\right)_\infty \left(e^{-2i\vartheta};q\right)_\infty}{\sin\vartheta} H_n(\cos\vartheta \mid q) = 2i\left(e^{2i\vartheta};q\right)_\infty \left(e^{-2i\vartheta};q\right)_\infty \frac{H_{n+1}(\cos\vartheta \mid q)}{\sqrt{q^{n+1}}}.$$

This is equivalent with

$$\delta_q \frac{\left(e^{2i\vartheta};q\right)_\infty \left(e^{-2i\vartheta};q\right)_\infty}{\sin\vartheta}\sum_{n\geq 0}\frac{H_n(\cos\vartheta \mid q)}{(q;q)_n}t^n = 2i\left(e^{2i\vartheta};q\right)_\infty \left(e^{-2i\vartheta};q\right)_\infty \frac{1}{\sqrt{q}}\sum_{n\geq 0}\frac{H_{n+1}(\cos\vartheta \mid q)}{(q;q)_n}\left(\frac{t}{\sqrt{q}}\right)^n$$

$$= 2i\left(e^{2i\vartheta};q\right)_\infty \left(e^{-2i\vartheta};q\right)_\infty \frac{1}{\sqrt{q}}\left(e^{i\vartheta} + e^{-i\vartheta} - \frac{t}{\sqrt{q}}\right)\sum_{n\geq 0}\frac{H_n(\cos\vartheta \mid q)}{(q;q)_n}\left(\frac{t}{\sqrt{q}}\right)^n$$

and follows from an easy calculation:

$$\delta_q 2i \frac{(e^{2i\vartheta};q)_\infty (e^{-2i\vartheta};q)_\infty}{e^{i\vartheta}-e^{-i\vartheta}} \frac{1}{(e^{i\vartheta}t;q)_\infty (e^{-i\vartheta}t;q)_\infty}$$

$$= 2i \frac{(e^{2i\vartheta};q)_\infty (e^{-2i\vartheta};q)_\infty}{\left(q^{-\frac{1}{2}}e^{i\vartheta}t;q\right)_\infty \left(q^{-\frac{1}{2}}e^{-i\vartheta}t;q\right)_\infty} \left( \frac{(1-q^{-1}e^{-2i\vartheta})(1-q^{-\frac{1}{2}}e^{i\vartheta}t)}{(1-e^{2i\vartheta})(q^{\frac{1}{2}}e^{i\vartheta}-q^{-\frac{1}{2}}e^{-i\vartheta})} - \frac{(1-q^{-1}e^{2i\vartheta})(1-q^{-\frac{1}{2}}e^{-i\vartheta}t)}{(1-e^{-2i\vartheta})(q^{-\frac{1}{2}}e^{i\vartheta}-q^{\frac{1}{2}}e^{-i\vartheta})} \right)$$

$$= \frac{2i}{\sqrt{q}}\left(e^{i\vartheta}+e^{-i\vartheta}-\frac{t}{\sqrt{q}}\right) \frac{(e^{2i\vartheta};q)_\infty (e^{-2i\vartheta};q)_\infty}{\left(q^{-\frac{1}{2}}e^{i\vartheta}t;q\right)_\infty \left(q^{-\frac{1}{2}}e^{-i\vartheta}t;q\right)_\infty}$$

If we write again $v(\vartheta;q) = \frac{(e^{2i\vartheta};q)_\infty (e^{-2i\vartheta};q)_\infty}{\sin \vartheta}$ then (6.23) implies

$$H_n(\cos\vartheta \mid q) = \frac{q-1}{2} q^{\frac{n-1}{2}} v(\vartheta;q)^{-1} \Delta_q v(\vartheta;q) H_{n-1}(\cos\vartheta \mid q).$$

By iterating this formula we get the following analogue of the Rodrigues formula (1.19):

$$H_n(\cos\vartheta \mid q) = \left(\frac{q-1}{2}\right)^n q^{\frac{n(n-1)}{4}} v(\vartheta;q)^{-1} \left(\Delta_q\right)^n v(\vartheta;q). \tag{6.24}$$